\documentclass[portuguese,a4paper,twoside,10pt]{article}
\usepackage[english]{babel}
\usepackage[T1]{fontenc}
\usepackage[latin1]{inputenc}
\usepackage{amssymb,amsmath,epsfig,amsfonts,euscript}
\usepackage{dsfont}
\usepackage[usenames,dvipsnames]{color}

\newtheorem{teo}{Theorem}[section]
\newtheorem{pro}[teo]{Proposition}

\newtheorem{ex}{Example}

\newcounter{note}[section]
\newcounter{example}[section]

\newcommand{\en}{\mathbb{N}}

\newcommand{\er}{\mathbb{R}}

\newcommand{\dst}{\displaystyle}
\newcommand{\dem}{{\bf Proof }}
\newcommand{\fdem}{$\square$}

\newcommand{\nid}{\noindent}
\newcommand{\nn}{\nonumber}
\newcommand{\mb}{\mathbf}

\newcommand{\titulo}[1]{\begin{center}\mbox{} \\ \noindent \textit{\textbf{\Large #1}}\\\vspace{0.5cm}\end{center}}

\renewcommand{\abstract}[1]{{\small \noindent \textbf{Abstract:} #1\\}}

\newcommand{\keywords}[1]{{\small \noindent \textbf{Keywords:} #1\\}}

\begin{document}

\begin{center}
\titulo{Extremes of multivariate ARMAX processes}
\end{center}

\vspace{0.5cm}

\textbf{Marta Ferreira} Center of Mathematics of Minho University/DMA, Braga, Portugal\\

\textbf{Helena Ferreira} Department of Mathematics, University of
Beira
Interior, Covilhã, Portugal\\

\abstract{We define a  new multivariate time series model by generalizing the ARMAX process in a multivariate way.
We give conditions on stationarity and analyze  local  dependence and domains of attraction. As a consequence of the obtained result, we derive a new method of construction of multivariate extreme value copulas. We characterize the extremal dependence by computing
the multivariate extremal index and bivariate upper tail dependence coefficients. An estimation procedure for the multivariate extremal index shall be presented. We also address the marginal estimation and propose a new estimator for the ARMAX autoregressive parameter. }


\section{Introduction}
Stationary time series presenting sudden large peaks are usually well modeled by heavy tailed noise ARMA. However, models with practical application but simpler treatment have been studied in literature as an alternative. Davis and Resnick (\cite{dav+res}, 1989) proposed the MARMA process which is analogous to the ARMA  by just replacing summation by the maximum operation:
$$
X_i=\phi_1X_{i-1}\vee ...\vee \phi_pX_{i-p}\vee Y_i\vee\theta_1Y_{i-1}\vee...\vee\theta_qY_{i-q},
$$
where $0\leq \phi_i,\theta_j\leq 1$, $1\leq i\leq p,1\leq j\leq q$ and the innovations $Y_n$, $n\geq 1$,  are independent with unit Fréchet distribution. A first order MARMA type version, i.e.,
\begin{eqnarray}\label{armax}
X_i=cX_{i-1}\vee  Y_i,\,\,0<c<1,
\end{eqnarray}
was analyzed in Alpuim (\cite{alp1}, 1989) by considering innovations $Y_n$, $n\geq 1$,  independent and equally  distributed, not necessarily unit Fréchet.
The model in (\ref{armax}), sometimes denoted in literature as ARMAX,  corresponds to the case $\alpha=0$ in the Haslett (\cite{has}, 1979) model
$$
X_i=\beta X_{i-1}\vee  (\alpha\beta X_{i-1}+ Y_i),
$$
$0<\beta<1, 0\leq \alpha\leq 1$, used to describe a solar thermal energy storage system and later developed in, e.g., Daley and Haslett (\cite{dal+has}1982) and Greenwood and Hooghiemstra (\cite{green+hoog}, 1988). Exponent versions of ARMAX, namely pARMAX and pRARMAX, were used in the modeling  of financial series (Ferreira and Canto e Castro, \cite{lccmf2} 2010). Further applications of ARMAX processes and their generalizations can be seen, for instance, in Lebedev (\cite{leb}, 2008) and references therein. Here we shall consider a multivariate formulation of ARMAX and extend some of the results in Alpuim (\cite{alp1}, 1989). More precisely, we will analyze conditions on stationarity (Section \ref{smodel}), local  dependence conditions (Section \ref{sdep}) and domains of attraction (Section \ref{sindext}). The relation between the max-attractors of the process and the innovations allow us to evidence a new construction method of copulas for multivariate extreme value distributions (MEV). In computing the multivariate extremal index, we shall find that it is possible to have clustering in all or only in some of the marginals, according to their domains of attraction. An estimation procedure for the multivariate extremal index will be also stated  (Section \ref{sindext}). In Section \ref{sdepindep} we will derive the lag-$r$ {tail dependence coefficient (TDC)} (Sibyua \cite{sib}, 1960; Joe \cite{joe}, 1997) and the lag-$r$ {tail independence coefficient} of Ledford and Tawn (\cite{led+tawn1,led+tawn1}, 1996/97) and we find different types of tail dependence. Some notes on the marginal parameters estimation shall be given at the end (Section \ref{sestim}). In particular, we present a new estimator for the ARMAX parameter $c$ which is strongly consistent and asymptotically normal.

\section{The multivariate model}\label{smodel}
Let $\{\mb{X}_n=(X_{n,1},...,X_{n,d})\}_{n\geq 1}$ be a d-variate sequence, such that
\begin{eqnarray}\label{mar}
X_{n,j}=c_jX_{n-1,j}\vee Y_{n,j},\,\,n\geq 1,\,\,j=1,...,d,\,\,0<c_j<1,
\end{eqnarray}
where $\mb{X}_0=(X_{0,1},...,X_{0,d})$, $\mb{Y}_n=(Y_{n,1},...,Y_{n,d})$, $n\geq 1$, are independent, $\mb{X}_0\sim F_0$ and $\mb{Y}_n\sim G$. $\{\mb{X}_n\}_{n\geq 1}$ thus corresponds to a $d$-variate formulation of an ARMAX process given in (\ref{armax}).

It is an immediate consequence of relation (\ref{mar}) that each marginal $\{X_{n,j}\}_{n\geq 1}$, $j\in D=\{1,...d\}$, of the sequence $\{\mb{X}_n\}_{n\geq 1}$ can be written in the form
\begin{eqnarray}\label{mar1}
X_{n,j}=c_j^nX_{0,j}\vee\bigvee_{i=1}^nc_j^{n-i} Y_{i,j},\,\,n\geq 1.
\end{eqnarray}
If $F_n$ denotes the distribution of $\mb{X}_n=(X_{n,1},...,X_{n,d})$, we also have
\begin{eqnarray}\label{fnmar}
F_{n}(x_1,...,x_d)=F_{n-1}\Big(\frac{x_1}{c_1},...,\frac{x_d}{c_d}\Big)G(x_1,...,x_d)
\end{eqnarray}
and
\begin{eqnarray}\label{fnmar1}
F_{n}(x_1,...,x_d)=F_{0}\Big(\frac{x_1}{c_1^n},...,\frac{x_d}{c_d^n}\Big)
\prod_{i=1}^nG\Big(\frac{x_1}{c_1^{n-i}},...,\frac{x_d}{c_d^{n-i}}\Big).
\end{eqnarray}

If there exists  $\lim_{n\to\infty}F_n(x_1,...,x_d)$ then, based on (\ref{fnmar}), we have
\begin{eqnarray}\label{fmar}
F(x_1,...,x_d)=F\Big(\frac{x_1}{c_1},...,\frac{x_d}{c_d}\Big)G(x_1,...,x_d).
\end{eqnarray}
and, by (\ref{fnmar1}),
\begin{eqnarray}\label{fmar1}
F(x_1,...,x_d)=\lim_{n\to\infty}
\prod_{i=1}^nG\Big(\frac{x_1}{c_1^{i}},...,\frac{x_d}{c_d^{i}}\Big)
=\prod_{i=1}^\infty G\Big(\frac{x_1}{c_1^{i}},...,\frac{x_d}{c_d^{i}}\Big).
\end{eqnarray}

\begin{pro} \label{pnd}
$\{\mb{X}_n\}_{n\geq 1}$ is a stationary sequence with common non-degenerate distribution if and only if there exists $(x_1,...,x_d)\in\er_+^d$ such that
\begin{eqnarray}\label{pnd1}
0<\sum_{i=1}^\infty -\log G\Big(\frac{x_1}{c_1^{i}},...,\frac{x_d}{c_d^{i}}\Big)<\infty.
\end{eqnarray}
In this case, the common distribution F of $\{\mb{X}_n\}_{n\geq 1}$ satisfies (\ref{fmar}).
\end{pro}
\dem By (\ref{fmar1}), $F$ is non degenerate if and only if there exists $\mb{x}=(x_1,...,x_d)$ such that $0<F(x_1,...,x_d)<1$, i.e., such that
\begin{eqnarray}\nn
0<\prod_{i=1}^\infty G\Big(\frac{x_1}{c_1^{i}},...,\frac{x_d}{c_d^{i}}\Big)<1.
\end{eqnarray}
The assertion in (\ref{pnd1}) is straightforward by taking logarithms and, if it holds for some $\mb{x}\in\er^d$, then $G\Big(\frac{x_1}{c_1^{i}},...,\frac{x_d}{c_d^{i}}\Big)\dst\mathop{\to}_{i\to\infty} 1$, and thus  $\mb{x}\in\er^d_+$. \fdem\\

As a consequence of (\ref{pnd1}), if any of the marginals $G_j$, $j\in D$, of $G$ has support with non positive right end-point then the corresponding marginal $F_{n,j}$ has degenerate limiting distribution and, therefore, $F$ is a $d$-dimensional degenerate distribution. Observe that (\ref{pnd1}) is satisfied by every multivariate distribution with positive dependence and marginals $G_j$ satisfying $0<\sum_{i=1}^{\infty}-\log G_j(x/c_j^i)<\infty$ for some $x>0$. This latter condition is satisfied, for instance, by the Generalized Pareto distribution (Alpuim, \cite{alp1} 1989).

For each $j\in D$, suppose that $F_j$ belongs to the max-domain of attraction of $H_j$, in short $F_j\in\mathcal{D}(H_j)$, i.e., there exists constants $\{a_{n,j}>0\}_{n\geq 1}$ and $\{b_{n,j}\}_{n\geq 1}$, such that
\begin{eqnarray}\nn
n(1-F(a_{n,j}x+b_{n,j}))\mathop{\to}_{n\to\infty} -\log H_j(x),
\end{eqnarray}
where $H_j$ may be a Gumbel, a Weibull or a Fréchet distribution, respectively, $\Lambda(x)=\exp(-e^{-x})$, $\Psi_{\alpha_j}(x)=e^{-(-x)^{\alpha_j}}$, $x\leq 0$, and $\Phi_{\alpha_j}(x)=e^{-x^{-\alpha_j}}$, $x>0$, for some $\alpha_j>0$. Therefore, a sequence of normalized levels $\{u_{n,j}^{(\tau_j)}\}_{n\geq 1}$ for $\{X_{n,j}\}_{n\geq 1}$, i.e., such that
\begin{eqnarray}\nn
n(1-F(u_{n,j}^{(\tau_j)}))\mathop{\to}_{n\to\infty} \tau_j\geq 0
\end{eqnarray}
can be written as
\begin{eqnarray}\nn
u_{n,j}^{(\tau_j)}=a_{n,j}H_j^{-1}(e^{-\tau_j})+b_{n,j},
\end{eqnarray}
with $H_j^{-1}(x)=\inf\{y:F(y)\geq x\}$ the generalized inverse of $H_j$. By applying the Khintchine's types theorem (see, e.g., Leadbetter \emph{et al.}, \cite{lead+} 1983), we arrive at the following property of the normalized levels for $\{X_{n,j}\}_{n\geq 1}$  that shall be used latter:
\begin{eqnarray}\label{tau*}
\frac{u_{n,j}^{(\tau_j)}}{c_j}=u_{n,j}^{(\tau_j^*)}, \,\textrm{ with }\,
\tau_j^*=\left\{
\begin{array}{ll}
0&,\, \textrm{if $H_j\in\{\Lambda, \Psi_{\alpha_j}\}$}\vspace{0.15cm}\\
\tau_j c_j^{\alpha_j}&,\, \textrm{if $H_j=\Phi_{\alpha_j}$}.
\end{array}
\right.
\end{eqnarray}
In the sequel we shall denote $\{\mb{u}_n^{(\boldsymbol{\tau})}=(u_{n,1}^{(\tau_1)},...,u_{n,d}^{(\tau_d)})\}_{n\geq 1}$ the sequence of normalized random vectors.

\section{Asymptotic independence and local dependence of $\{\mb{X}_n\}_{n\geq 1}$}\label{sdep}

As showed in Alpuim (\cite{alp1}, 1989) for the univariate case, we shall prove that the strong-mixing condition also holds for the multivariate sequence, i.e., for any $A\in \mathcal{B}(\mb{X}_1,...,\mb{X}_p)$ and $B\in \mathcal{B}(\mb{X}_{p+s+1},\mb{X}_{p+s+2},...)$,
\begin{eqnarray}\nn
|P(A\cap B)-P(A)P(B)|\leq \alpha_s
\end{eqnarray}
with $\alpha_s\dst\mathop{\to}_{s\to\infty} 0$, where $\mathcal{B}(\cdot)$ denotes the $\sigma$-field generated by the indicated random vectors.

In what follows, all operations and inequalities 
 between vectors are understood to be componentwise.

\begin{pro} \label{psmix}
$\{\mb{X}_n\}_{n\geq 1}$  satisfies the strong-mixing condition.
\end{pro}
\dem Consider $A\in \mathcal{B}(\mb{X}_1,...,\mb{X}_p)$ and $B\in \mathcal{B}(\mb{X}_{p+s+1},\mb{X}_{p+s+2},...)$ and let
\begin{eqnarray}\nn
C_s=\{\mb{Y}_{p+1}\not\geq \mb{X}_{p},...,\mb{Y}_{p+s+1}\not\geq \mb{X}_{p+s}\}.
\end{eqnarray}

We have
\begin{eqnarray}\nn
\begin{array}{rl}
&|P(A\cap B)-P(A)P(B)|\vspace{0.35cm}\\
=&|P(A\cap B\cap C_s)+P(A\cap B\cap \overline{C}_s)-P(A)P(B\cap C_s)-P(A)P(B\cap \overline{C}_s)|\vspace{0.35cm}\\
\leq &|P(A\cap B\cap C_s)-P(A)P(B\cap C_s)|+|P(A\cap B\cap \overline{C}_s)-P(A)P(B\cap \overline{C}_s)| \end{array}
\end{eqnarray}

Observe that, for the first term,
\begin{eqnarray}\nn
\begin{array}{rl}
&|P(A\cap B\cap C_s)-P(A)P(B\cap C_s)|\vspace{0.35cm}\\
=& |P(C_s)P(B|C_s)P(A|B\cap C_s)-P(A)P(C_s)P(B|C_s)|\vspace{0.35cm}\\
\leq & P(C_s)P(B|C_s)|P(A|B\cap C_s)-P(B|C_s)|
\vspace{0.35cm}\\
\leq &P(C_s).
\end{array}
\end{eqnarray}
On the other hand, since
\begin{eqnarray}\nn
\mb{X}_{p+s+1}=\mb{c}^{p+s+1-k}\mb{X}_k\vee \bigvee_{i=1}^{p+s+1-k}\mb{c}^{p+s+1-k-i}\mb{Y}_{k+i},\,\,k=p,...,p+s,
\end{eqnarray}
we can write
\begin{eqnarray}\nn
A\cap B\cap \overline{C}_s=A\cap B'\cap \overline{C}_s \textrm{ and }  B\cap \overline{C}_s= B'\cap \overline{C}_s
\end{eqnarray}
where $B'\in \mathcal{B}(\mb{Y}_{p+1},\mb{Y}_{p+2},...)$. Thus being, for the second term, we have
\begin{eqnarray}\nn
\begin{array}{rl}
&|P(A\cap B\cap \overline{C}_s)-P(A)P(B\cap \overline{C}_s)|= |P(A\cap B'\cap \overline{C}_s)-P(A)P(B'\cap \overline{C}_s)|\vspace{0.35cm}\\
\leq & |P(A\cap B\cap \overline{C}_s)-P(A\cap B')|+|P(A\cap B')-P(A)P(B\cap \overline{C}_s)|
\vspace{0.35cm}\\
=& P(A\cap B')P({C}_s|A\cap B')+|P(A)P(B')-P(A)P(B\cap \overline{C}_s)|
\vspace{0.35cm}\\
\leq &2P(C_s).
\end{array}
\end{eqnarray}
Now we just need to prove that $P(C_s)\dst\mathop{\to}_{s\to\infty}0$. Observe that
\begin{eqnarray}\nn
\begin{array}{rl}
P(C_s)=&1-P(\bigcup_{i=1}^{s+1}\{\mb{Y}_{p+i}\geq \mb{X}_{p+i-1}\})\leq 1-P(\bigcup_{i=1}^{s+1}\{\mb{Y}_{p+i}\geq \mb{c}^{i-1}\mb{X}_{p}\})\vspace{0.35cm}\\
\leq &1-P(\mb{Y}_{p+s+1}\geq \mb{c}^{s}\mb{X}_{p})
=1-\int_{\er^d}H\big(\frac{y_1}{c_1^s},...,\frac{y_d}{c_d^s}\big)dG(y_1,...,y_d)
\dst\mathop{\to}_{s\to\infty}0.\,\,\square
\end{array}
\end{eqnarray}

Therefore, $\{\mb{X}_n\}_{n\geq 1}$ satisfies condition  $D(\mathbf{u}_n,\alpha_{l_n})$, for any sequence of real vectors $\{\mb{u}_n\}_{n\geq 1}$ and for any sequence $\{l_n\}_{n\geq 1}$ such that $l_n\to\infty$,  corresponding to the multivariate version of Leadbetter's D-condition of local dependence (see, e.g., Leadbetter \emph{et al.} \cite{lead+} 1983).

Now we shall see that $\{\mb{X}_n\}_{n\geq 1}$ also satisfies the multivariate version of $D''$ condition of Leadbetter and Nandagopalan (\cite{lead+nand}, 1989). For a given sequence of real vectors $\{\mb{u}_n\}_{n\geq 1}$, we say that condition $D''(\mb{u}_n)$ holds if $D(\mathbf{u}_n,\alpha_{l_n})$ also holds and
\begin{eqnarray}\nn
\begin{array}{rl}
n\sum_{i=2}^{[n/k_n]}P(\mb{X}_1\not\leq \mb{u}_n,\mb{X}_i\leq \mb{u}_n,\mb{X}_{i+1}\not\leq \mb{u}_n)\to 0
\end{array}
\end{eqnarray}
for some sequence $\{k_n\}_{n\geq 1}$ such that, as $n\to\infty$,
\begin{eqnarray}\nn
\begin{array}{rl}
k_n\to\infty, \,\,\frac{k_nl_n}{n}\to 0 \,\,\textrm{ and }\,\, k_n\alpha_{l_n}\to 0.
\end{array}
\end{eqnarray}

\begin{pro}
$\{\mb{X}_n\}_{n\geq 1}$ satisfies condition $D''(\mb{u}_n^{(\boldsymbol{\tau})})$.
\end{pro}
\dem Observe that
\begin{eqnarray}\nn
\begin{array}{rl}
&n\sum_{i=2}^{[n/k_n]}P(\mb{X}_1\not\leq \mb{u}_n^{(\boldsymbol{\tau})},\mb{X}_i\leq \mb{u}_n^{(\boldsymbol{\tau})},\mb{X}_{i+1}\not\leq \mb{u}_n^{(\boldsymbol{\tau})})\vspace{0.35cm}\\
\leq & \sum_{j=1}^dn\sum_{i=2}^{[n/k_n]}P({X}_{1,j}> {u}_{n,j}^{(\tau_j)}, {X}_{i,j}\leq {u}_{n,j}^{(\tau_j)}<{X}_{i+1,j} )\vspace{0.35cm}\\
&+\sum_{1\leq s,s'\leq d}n\sum_{i=2}^{[n/k_n]}P({X}_{1,s}> {u}_{n,s}^{(\tau_s)},\mb{X}_i\leq \mb{u}_n^{(\boldsymbol{\tau})},{X}_{i+1,s'}>{u}_{n,s'}^{(\tau_{s'})}).
\end{array}
\end{eqnarray}
Since each marginal sequence $\{X_{n,j}\}_{n\geq 1}$ satisfies condition $D''(u_{n,j}^{(\tau_j)})$ (Canto e Castro, \cite{lcc} 1992), the first term above has null limit, as $n\to\infty$. 
The second term above is upper bounded by, successively,
\begin{eqnarray}\nn
\begin{array}{rl}
&\sum_{1\leq s,s'\leq d}n\sum_{i=2}^{[n/k_n]}P({X}_{1,s}> {u}_{n,s}^{(\tau_s)},{X}_{i,s'}\leq {u}_{n,s'}^{(\tau_{s'})}<c_{s'}X_{i,s'}\vee Y_{i+1,s'})\vspace{0.35cm}\\
= &\sum_{1\leq s,s'\leq d}n\sum_{i=2}^{[n/k_n]}P({X}_{1,s}> {u}_{n,s}^{(\tau_s)},{X}_{i,s'}\leq {u}_{n,s'}^{(\tau_{s'})}< Y_{i+1,s'}) \vspace{0.35cm}\\
\leq &\sum_{1\leq s,s'\leq d}n\sum_{i=2}^{[n/k_n]}P({X}_{1,s}> {u}_{n,s}^{(\tau_s)})P( Y_{i+1,s'}> {u}_{n,s'}^{(\tau_{s'})})\vspace{0.35cm}\\
= &\sum_{1\leq s,s'\leq d}n\big[\frac{n}{k_n}\big](1-F_s({u}_{n,s}^{(\tau_s)}))
\Big(1-\frac{F_{s'}({u}_{n,s'}^{(\tau_{s'})})}{F_{s'}({u}_{n,s'}^{(\tau_{s'})}/c_{s'})}\Big)
\vspace{0.35cm}\\
\leq &\frac{1}{k_n}\sum_{1\leq s,s'\leq d}n(1-F_s({u}_{n,s}^{(\tau_s)}))\Big(n(1-F_{s'}({u}_{n,s'}^{(\tau_{s'})}))-
n(1-F_{s'}({u}_{n,s'}^{(\tau_{s'})}/c_{s'}))\Big)\frac{1}{F_{s'}({u}_{n,s'}^{(\tau_{s'})}/c_{s'})},
\end{array}
\end{eqnarray}
which also converges to zero for any sequence $k_n\to\infty$, since by (\ref{tau*}) we have $n(1-F_{s'}({u}_{n,s'}^{(\tau_{s'})}/c_{s'}))\dst\mathop{\to}_{} \tau^*_{s'}\geq 0$, as $n\to\infty$. \fdem\\

\section{The multivariate extremal index and the domain of attraction of $\{\mb{X}_n\}_{n\geq 1}$}\label{sindext}

\begin{figure}
\begin{center}
\includegraphics[width=12cm,height=9cm]{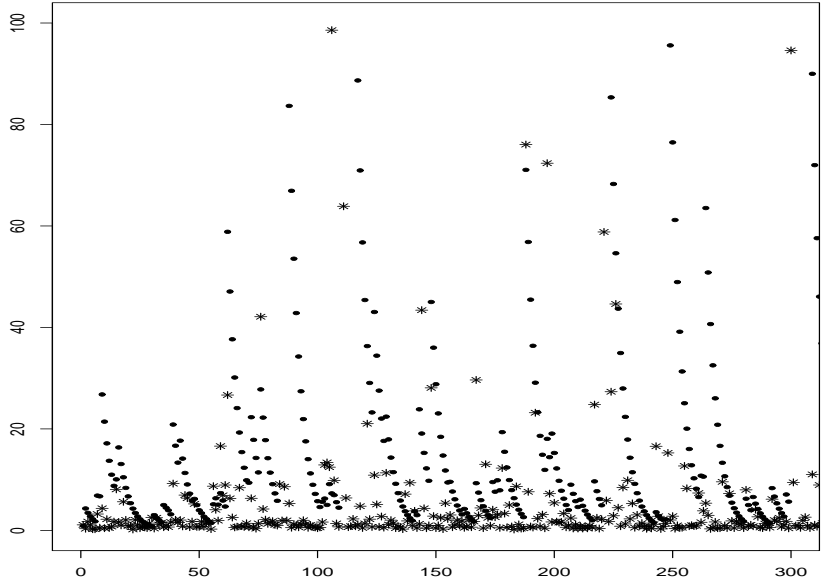}
\end{center}
\caption{Marginal sample paths of a bivariate ARMAX process with unit Fréchet innovations: the symbol "$\ast$" corresponds to c=$0.1$ where the high values tend to occur almost singly and symbol "$\bullet$" corresponds to c=$0.8$ with visible high values clustering. \label{fig3}}
\end{figure}

A phenomenon also noticed in real data is that extreme events often tend to occur in clusters. The measure that is used to capture the clustered extremal dependence is the \emph{extremal index} (Leadbetter \emph{et al.} \cite{lead+} 1983). More precisely, the  extremal index can be interpreted as the reciprocal of the limiting mean cluster size. A unit extremal index means no serial clustering and is a form of asymptotic independence of extremes. Figure \ref{fig3} illustrates both clustering (symbol "$\bullet$") and asymptotic independence (symbol "$\ast$") at high levels.\\

The results of the previous section will allow us  to compute the multivariate extremal index of $\{\mb{X}_n\}_{n\geq 1}$ (Nandagopalan \cite{nand} 1990). More precisely, if for all $\boldsymbol{\tau}\in\er_+^d$ there exists normalized levels $\{\mb{u}_n^{(\boldsymbol{\tau})}=(u_{n,1}^{(\tau_1)},...,u_{n,d}^{(\tau_d)})\}_{n\geq 1}$, such that the sequence $\{nP(\mb{X}_1\not\leq \mb{u}_n^{(\boldsymbol{\tau})})\}_{n\geq 1}$ is convergent and $D''(\mb{u}_n^{(\boldsymbol{\tau})})$ holds, then $\{\mb{X}_n\}_{n\geq 1}$ has multivariate extremal index if and only if, $\forall \boldsymbol{\tau}\in \er_+^d$, sequence $\{nP(\mb{X}_1\leq \mb{u}_n^{(\boldsymbol{\tau})},\mb{X}_2\not\leq \mb{u}_n^{(\boldsymbol{\tau})})\}_{n\geq 1}$ converges too. In this case,
\begin{eqnarray}\label{indextm}
\begin{array}{rl}
\theta(\tau_1,...,\tau_d)=\lim_{n\to\infty}\frac{P(\mb{X}_1\leq \mb{u}_n^{(\boldsymbol{\tau})},\mb{X}_2\not\leq \mb{u}_n^{(\boldsymbol{\tau})})}{P(\mb{X}_1\not\leq \mb{u}_n^{(\boldsymbol{\tau})})},\,\boldsymbol{\tau}\in\er_+^d,
\end{array}
\end{eqnarray}
(Ferreira, \cite{hf94} 1994).
In the sequel we shall use the copula function with notation
$$
C_F(u_1,..,u_d)=F(F_1^{-1}(u_1),...,F_d^{-1}(u_d)),\,(u_1,...,u_d)\in[0,1]^d.
$$

\begin{pro}
If $F\in \mathcal{D}(H)$ then $\{\mb{X}_n\}_{n\geq 1}$ has  multivariate extremal index
\begin{eqnarray}\nn
\begin{array}{rl}
\dst\theta(\tau_1,...,\tau_d)=1-\frac{\log C_{H_I}(e^{-\tau_jc_j^{\alpha_j}},j\in I)}{\log C_{H}(e^{-\tau_1},...,e^{-\tau_d})},
\end{array}
\end{eqnarray}
where $I$ is the set of indexes in $D$ for which $H_j(x)=\Phi_{\alpha_j}(x)=e^{-x^{-\alpha_j}}$, $x>0$. Moreover,
\begin{eqnarray}\label{pindext}
\theta_j=\left\{
\begin{array}{ll}
1&,\, \textrm{if $H_j\in\{\Lambda, \Psi_{\alpha_j}\}$}\vspace{0.15cm}\\
1-c_j^{\alpha_j}&,\, \textrm{if $H_j=\Phi_{\alpha_j}$}
\end{array}
\right.
\end{eqnarray}
is the extremal index of $\{X_{n,j}\}_{n\geq 1}$, $j=1,...,d$.
\end{pro}
\dem By hypothesis, $F\in \mathcal{D}(H)$, i.e., $F_j\in \mathcal{D}(H_j,\{a_{n,j}>0\},\{b_{n,j}\})$, with $H_j$ of the extremal type $\Lambda$, $\Psi_{\alpha_j}$ or $\Phi_{\alpha_j}$, and $C_F^n(u_1^{1/n},..,u_d^{1/n})\dst\mathop{\to}_{n\to\infty}C_H(u_1,...,u_d)$, $(u_1,...,u_d)\in[0,1]^d$. Thus we guarantee the existence of normalized levels $\mb{u}_n^{(\boldsymbol{\tau})}=(u_{n,1}^{(\tau_1)},...,u_{n,d}^{(\tau_d)})$, for which condition $D''(\mb{u}_n^{(\boldsymbol{\tau})})$ holds.

Moreover, for $u_{n,j}^{(\tau_j)}=a_{n,j}H_j^{-1}(e^{-\tau_j})+b_{n,j}$, $j=1,...,d$, we have
\begin{eqnarray}\nn
\begin{array}{rl}
nP(\mb{X}_1\not\leq \mb{u}_n^{(\boldsymbol{\tau})})=&n(1-F(a_{n,1}H_1^{-1}(e^{-\tau_1})+b_{n,1},...,
a_{n,d}H_d^{-1}(e^{-\tau_d})+b_{n,d}))\vspace{0.35cm}\\
 \dst\mathop{\to}_{n\to\infty}&-\log H(H_1^{-1}(e^{-\tau_1}),...,H_d^{-1}(e^{-\tau_d})).
\end{array}
\end{eqnarray}
On the other hand,
\begin{eqnarray}\nn
\begin{array}{rl}
&nP(\mb{X}_1\leq \mb{u}_n^{(\boldsymbol{\tau})},\mb{X}_2\not\leq \mb{u}_n^{(\boldsymbol{\tau})})
=
nP(\mb{X}_1\leq \mb{u}_n^{(\boldsymbol{\tau})},\mb{c}\mb{X}_1\vee \mb{Y}_2\not\leq \mb{u}_n^{(\boldsymbol{\tau})})\vspace{0.35cm}\\
=&
nP(\mb{X}_1\leq \mb{u}_n^{(\boldsymbol{\tau})})P(\mb{Y}_2\not\leq \mb{u}_n^{(\boldsymbol{\tau})})=P(\mb{X}_1\leq \mb{u}_n^{(\boldsymbol{\tau})})n(1-G(u_{n,1}^{(\tau_1)},...,u_{n,d}^{(\tau_d)}))
\vspace{0.35cm}\\
=&\frac{P(\mb{X}_1\leq \mb{u}_n^{(\boldsymbol{\tau})})}{P(\mb{X}_1\leq \mb{u}_n^{(\boldsymbol{\tau})}/\mb{c})}\Big(n(1-F(u_{n,1}^{(\tau_1)},...,u_{n,d}^{(\tau_d)}))
-n\big(1-F\big(\frac{u_{n,1}^{(\tau_1)}}{c_1},...,\frac{u_{n,d}^{(\tau_d)}}{c_d}\big)\big)\Big)
\vspace{0.35cm}\\
\dst\mathop{\to}_{n\to\infty}& -\log H(H_1^{-1}(e^{-\tau_1}),...,H_d^{-1}(e^{-\tau_d}))+\log H(H_1^{-1}(e^{-\tau_1^*}),...,H_d^{-1}(e^{-\tau_d^*}))\vspace{0.35cm}\\
=& -\log H(H_1^{-1}(e^{-\tau_1}),...,H_d^{-1}(e^{-\tau_d}))+\log H_I(H_1^{-1}(e^{-\tau_1^*}),...,H_d^{-1}(e^{-\tau_d^*}))_I,
\end{array}
\end{eqnarray}
where $I$ is the set of indexes in $D$ for which $\tau_j^*$ given in (\ref{tau*}) are positive, i.e., for which $H_j(x)=\Phi_{\alpha_j}(x)=e^{-x^{-\alpha_j}}$, $x>0$, and $H_I$ denotes the marginal distribution of $H$ corresponding to those indexes. Therefore, applying (\ref{indextm}), we have
\begin{eqnarray}\nn
\begin{array}{rl}
\theta(\tau_1,...,\tau_d)=&\dst 1-\frac{\log {H_I}(H_1^{-1}(e^{-\tau_1^*}),...,H_d^{-1}(e^{-\tau_d^*}))_I}{\log {H}(H_1^{-1}(e^{-\tau_1}),...,H_d^{-1}(e^{-\tau_d}))}\vspace{0.35cm}\\
=&\dst1-\frac{\log C_{H_I}(e^{-\tau_1^*},...,e^{-\tau_d^*})_I}{\log C_{H}(e^{-\tau_1},...,e^{-\tau_d})}.
\end{array}
\end{eqnarray}
Observe that if $I=\emptyset$ then $\theta(\tau_1,...,\tau_d)=1$, $\forall \boldsymbol{\tau}$, and if $I\not =\emptyset$, we have
\begin{eqnarray}\nn
\theta_j=\left\{
\begin{array}{ll}
1&,\, \textrm{if $j\in D-I$ }\vspace{0.15cm}\\
1-\frac{\tau_j^*}{\tau_j}&,\, \textrm{if $j\in I$}
\end{array}
\right.
\end{eqnarray}
leading to the assertion (\ref{pindext}), which corresponds to the univariate marginal extremal index already derived in Alpuim (\cite{alp1}, 1989). \fdem\\

The expression obtained for the multivariate extremal index function has the advantage of, once known/estimated the constants $c_j$ and the marginal domains of attraction, we  are only dependent on the attractor copula of $\{\mb{\widehat{X}}_n\}_{n\geq 1}$ corresponding to the i.i.d.\,sequence with the same distribution $F$.

Since we have $C_F^n(u_1^{1/n},...,u_d^{1/n})\dst\mathop{\to}_{n\to\infty} C_H(u_1,...,u_d)$ uniformly in $[0,1]^d$ and $C_H$ is continuous, we can replace the discrete variable $n$ by a continuous variable $t$ and equivalently state $t(1-C_F(1-x_1/t,\hdots,1-x_d/t))\dst\mathop{\to}_{t\to\infty} \log C_H(e^{-x_1},...,e^{-x_d})$, $\mb{x}\in[0,\infty)^d$. If we rewrite the result of the previous proposition as
\begin{eqnarray}\nn
\begin{array}{rl}
\theta(\tau_1,...,\tau_d)=&\dst 1-\lim_{t\to\infty}\frac{t(1-C_{F_I}(1-\tau_1c_1^{\alpha_1}/t,\hdots,  1-\tau_dc_d^{\alpha_d}/t)_I) }{t(1-C_{F}(1-\tau_1/t,\hdots,  1-\tau_d/t))}\vspace{0.35cm}\\
=&\dst 1-\lim_{t\to\infty}\frac{t P(\bigcup_{j\in I}\{F_{X_{1,j}}(X_{1,j})>1-\tau_1c_1^{\alpha_1}/t\}) }{t P(\bigcup_{j=1}^d\{F_{X_{1,j}}(X_{1,j})>1-\tau_1/t\})},
\end{array}
\end{eqnarray}
we can then estimate the multivariate extremal index through tail dependence functions estimators concerning  $F_I$ and $F$. For this issue see, e.g., Schmidt and Stadtm\"{u}ller (\cite{schm+stadt}, 2006), Einmahl \emph{et al.} (\cite{einmahl+12}, 2012) and references therein.

\begin{ex}
Consider $F$ with $F_1,F_2\in\mathcal{D}(\Lambda)$ and $F_j\in\mathcal{D}(\Phi_1)$, $j=3,...,d$. If $C_F(u_1,...,u_d)=\exp(-(\sum_{j=1}^d(-\log u_j)^\gamma)^{1/\gamma})$, $\gamma \geq 1$, then $F\in\mathcal{D}(H)$, with $C_H=C_F$, $H_1=H_2=\Lambda$ and $H_j=\Phi_1$, $j=3,...,d$. Therefore, we have
\begin{eqnarray}\nn
\begin{array}{rl}
\dst\theta(\tau_1,...,\tau_d)=1-\frac{(\sum_{j=3}^d(\tau_jc_j)^\gamma)^{1/\gamma}}
{(\sum_{j=1}^d\tau_j^\gamma)^{1/\gamma}},\,\,(\tau_1,...,\tau_d)\in\er_+^d,
\end{array}
\end{eqnarray}
$\theta_1=\theta_2=1$ and $\theta_j=1-c_j$, $j=3,...,d$.
\end{ex}

\begin{ex}
Consider $F$ with $F_j\in\mathcal{D}(\Phi_1)$, $j=1,...,d$ and $C_F(u_1,...,u_d)=\bigwedge_{j=1}^d u_j$. Then $F\in\mathcal{D}(H)$, with $C_H=C_F$, $H_j=\Phi_1$, $j=1,...,d$. Therefore, we have
\begin{eqnarray}\nn
\begin{array}{rl}
\dst\theta(\tau_1,...,\tau_d)=1-\frac{\bigvee_{j=1}^d\tau_jc_j}
{\bigvee_{j=1}^d\tau_j},\,\,(\tau_1,...,\tau_d)\in\er_+^d.
\end{array}
\end{eqnarray}
\end{ex}

The next result relates the domain of attraction of $F$ with the one of $G$.

\begin{pro}\label{pcoprelFG}
If $F\in\mathcal{D}(H)$ then $G\in \mathcal{D}(V)$ with $V_j=H_j^{\theta_j}$ and $\theta_j$ given in (\ref{pindext}), $j\in D$,
and
\begin{eqnarray}\label{pcoprelFG1}
C_V(u_1,...,u_d)=\frac{C_H(u_1^{1/\theta_1},...,u_d^{1/\theta_d})}
{C_H(u_1^{1/\theta_1-1},...,u_d^{1/\theta_d-1})}
\end{eqnarray}
\end{pro}
\dem By hypothesis, $F_j\in \mathcal{D}(H_j,\{a_{n,j}>0\},\{b_{n,j})\})$, $j\in D$, i.e.,
$F_j^n(a_{n,j}x_j+b_{n,j})\dst\mathop{\to}_{n\to\infty} H_j(x_j)$
and $C_F^n(u_1^{1/n},..,u_d^{1/n})\dst\mathop{\to}_{n\to\infty}C_H(u_1,...,u_d)$, $(u_1,...,u_d)\in[0,1]^d$. In addition,
\begin{eqnarray}\nn
F_j^n\big(\frac{a_{n,j}x_j+b_{n,j}}{c_j}\big)\dst\mathop{\to}_{n\to\infty}
\left\{
\begin{array}{ll}
1&,\,\, \textrm{if $H_j\in\{\Lambda, \Psi_{\alpha_j}\}$}\vspace{0.15cm}\\
H_j^{c_j^{\alpha_j}}(x_j)&,\, \textrm{if $H_j=\Phi_{\alpha_j}$}.
\end{array}
\right.
\end{eqnarray}
From the stationarity relation in (\ref{fmar}), we have
\begin{eqnarray}\nn
F_j^n(a_{n,j}x_j+b_{n,j})=F_j^n\big(\frac{a_{n,j}x_j+b_{n,j}}{c_j}\big)G_j^n(a_{n,j}x_j+b_{n,j}).
\end{eqnarray}
Therefore,
\begin{eqnarray}\nn
G_j^n(a_{n,j}x_j+b_{n,j})\dst\mathop{\to}_{n\to\infty}
\left\{
\begin{array}{ll}
H_j(x_j)&,\,\, \textrm{if $H_j\in\{\Lambda, \Psi_{\alpha_j}\}$}\vspace{0.15cm}\\
H_j^{1-c_j^{\alpha_j}}(x_j)&,\, \textrm{if $H_j=\Phi_{\alpha_j}$},
\end{array}
\right.
\end{eqnarray}
and thus $G_j\in\mathcal{D}(H_j^{\theta_j})$, $j\in D$, with $\theta_j$ given in (\ref{pindext}).\\

\nid Now we look at the copula of $G$. We have
\begin{eqnarray}\nn
F^n(a_{n,1}x_1+b_{n,1},...,a_{n,d}x_d+b_{n,d})\dst\mathop{\to}_{n\to\infty}
H(x_1,...,x_d)=C_H(H_1(x_1),...,H_d(x_d))
\end{eqnarray}
and
\begin{eqnarray}\nn
\begin{array}{rl}
F^n\Big(\frac{a_{n,1}x_1+b_{n,1}}{c_1},...,
\frac{a_{n,d}x_1+b_{n,d}}{c_d}\Big)=&
C_F^n\Big(\Big(F_1^n\Big(\frac{a_{n,1}x_1+b_{n,1}}{c_1}\Big)\Big)^{1/n},...,
\Big(F_d^n\Big(\frac{a_{n,d}x_1+b_{n,d}}{c_d}\Big)\Big)^{1/n}\Big)\vspace{0.35cm}\\
\dst\mathop{\to}_{n\to\infty}&C_{H_I}\Big(H_j^{c_j^{\alpha_j}}(x_j),j\in I\Big)=C_{H_I}\Big(H_j\big(\frac{x_j}{c_j}\big),j\in I\Big)\vspace{0.35cm}\\
=&{H_I}\big(\frac{x_j}{c_j},j\in I\big)
\end{array}
\end{eqnarray}
Again from the relation between $F$ and $G$ in (\ref{fmar}), we obtain
\begin{eqnarray}\nn
G^n(a_{n,1}x_1+b_{n,1},...,a_{n,d}x_d+b_{n,d})\dst\mathop{\to}_{n\to\infty}
\frac{C_H(H_1(x_1),...,H_d(x_d))}{C_{H_I}\Big(H_j^{c_j^{\alpha_j}}(x_j),j\in I\Big)}.
\end{eqnarray}
Thus we can say that $G\in\mathcal{D}(V)$, where $V_j=H_j^{\theta_j}$ and
\begin{eqnarray}\nn
H(x_1,...,x_d)={H_I}\Big(\frac{x_j}{c_j},j\in I\Big)V(x_1,...,x_d),
\end{eqnarray}
or equivalently,
\begin{eqnarray}\nn
C_H(u_1,...,u_d)=C_{H}\big(u_1^{1-\theta_1},...,u_d^{1-\theta_d}\big)
C_V\big(u_1^{\theta_1},...,u_d^{\theta_d}\big),
\end{eqnarray}
given $\theta_j$, $j\in D$, stated in (\ref{pindext}). \fdem\\

Observe that if $\theta_j=\theta$, $j=1,\hdots,d$ or if $C_H$ is the product copula then $C_V=C_H$. However, in general,  relation (\ref{pcoprelFG1}) adds one more method to the existing ones of copulas construction (Joe \cite{joe} 1997, Liebschen \cite{lieb} 2008). Any MEV copula can appear in the limiting distribution $H$ generated from this model, since $H\in\mathcal{D}(H)$. Thus considering a MEV copula $C$ and constants $\theta_j\in]0,1]$, $j=1,\hdots,d$, we can derive new MEV copulas by applying the ratio rule in (\ref{pcoprelFG1}) one or more times. We shall illustrate the procedure by considering that $C$ is a Gumbel copula. More precisely, if $C(u_1,...,u_d)=\exp(-(\sum_{j=1}^d(-\log u_j)^{\gamma})^{1/\gamma})$, $\gamma\geq 1$, then
$$
C^*(u_1,...,u_d)=\frac{\exp(-(\sum_{j=1}^d(-\frac{1}{\theta_j}\log u_j)^{\gamma})^{1/\gamma})}{\exp(-(\sum_{j=1}^d(-(\frac{1}{\theta_j}-1)\log u_j)^{\gamma})^{1/\gamma})}
$$
is a MEV copula. The extremal coefficients, and thus the tail behavior, of this new copula present a greater variability of values when compared with the respective ones of the Gumbel. Indeed, for any $J\subset D$, we have $C_J(u,\hdots,u)=u^{\epsilon_J^C}$ with $\epsilon_J^C=|J|^{1/\gamma}$, and $C^*_J(u,\hdots,u)=u^{\epsilon_J^{C^*}}$ with $\epsilon_J^{C^*}=(\sum_{j\in J}(\frac{1}{\theta_j})^{\gamma})^{1/\gamma}-(\sum_{j\in J}(\frac{1}{\theta_j}-1)^{\gamma})^{1/\gamma}$.

\section{Coefficients of tail dependence and tail independence}\label{sdepindep}

Loosely speaking, tail
dependence describes the limiting proportion of exceedances by one margin of a certain high
threshold given that the other margin has already exceeded that threshold.  The most used definition of tail dependence, provided in the monograph of Joe (\cite{joe}, 1997), is the \emph{tail dependence coefficient (TDC)}:
\begin{eqnarray}\label{lambda}
\lambda=\dst\lim_{t\downarrow 0}P(F_Y(Y)>1-t|F_X(X)>1-t).
\end{eqnarray}
We say that the random pair $(X,Y)$ is tail dependent whenever $\lambda>0$ and tail independent if $\lambda=0$.

In the tail independent case, Ledford and Tawn (\cite{led+tawn1,led+tawn2} 1996/1997) proposed to model the
null limit in (\ref{lambda}) by introducing a coefficient ($\eta$) to rule the decay rate of the joint bivariate survival function:
\begin{eqnarray}\nn
P(F_Y(Y)>1-t|F_X(X)>1-t)\sim L(t)t^{1/\eta-1}, \textrm{ as $t\downarrow 0$,}
\end{eqnarray}
where $L$ is a slowly varying function at $0$, i.e. $L(tx)/L(t)\to 1$ as $t\downarrow 0$, for any
fixed $x > 0$  and $\eta\in (0,1]$ is a constant. Coefficient $\eta$ measures the degree of tail independence between r.v.'s $X$ and $Y$. Observe that tail dependence occurs if $\eta=1$ and $L(t)\not \to 0$, as $t\downarrow 0$, and tail independence otherwise. The r.v.'s $X$ and $Y$ are called positively associated when
$1/2 < \eta < 1$, nearly independent when $\eta= 1/2$ and negatively associated when $0 < \eta < 1/2$.

Both concepts can be naturally extended to a lag-$r$ ($r\in\en_0$) formulation of a stationary $d$-dimensional sequence, $\{\mb{X}_n=(X_{n,1},...,X_{n,d})\}_{n\geq 1}$. More precisely, the lag-$r$ TDC as
\begin{eqnarray}\nn
\lambda_{jj'}^{(r)}(\mb{X})=\dst\lim_{t\downarrow 0}P(F_{j'}(X_{1+r,j'}) > 1-t| F_j(X_{1,j}) > 1-t)
\end{eqnarray}
and the lag-$r$ ($r\in\en_0$) Ledford and Tawn  coefficient $\eta_{jj'}^{(r)}(\mb{X})$  defined by
\begin{eqnarray}\nn
P(F_{j'}(X_{1+r,j'}) > 1-t| F_j(X_{1,j}) > 1-t)\sim L(t)t^{1/\eta_{jj'}^{(r)}(\mb{X})-1}, \textrm{ as $t\downarrow 0$},
\end{eqnarray}
where $L$ is a slowly varying function at $0$. We denote $\lambda_{jj'}(\mb{X})\equiv \lambda_{jj'}^{(0)}(\mb{X})$ as the TDC between the $j$th and the $j'$th components, $\lambda_{j}^{(r)}(\mb{X})\equiv \lambda_{jj}^{(r)}(\mb{X})$ is the lag-$r$ TDC within the $j$th sequence and $\lambda_{jj'}^{(r)}(\mb{X})$
is the lag-$r$ cross-sectional TDC between the $j$th and the $j'$th sequences. An analogous description holds for the Ledford and Tawn coefficients, respectively, $\eta_{jj'}(\mb{X})$, $\eta_{j}^{(r)}(\mb{X})$ and $\eta_{jj'}^{(r)}(\mb{X})$.\\

In the sequel we shall denote $w_{j't}=F_{j'}^{-1}(1-t)$.

\begin{pro}\label{plambda1}
If $\{\mb{X}_n\}_{n\geq 1}$ has stationary distribution then
\begin{eqnarray}\nn
\begin{array}{rl}
\lambda_{jj'}^{(r)}(\mb{X})=2-\dst\lim_{t\downarrow 0}\frac{1}{t}\bigg(1-C_{jj'}( 1-t,F_{j'}(c_{j'}^{-r}w_{j't}))\frac{1-t}{F_{j'}(c_{j'}^{-r}w_{j't})}\bigg).
\end{array}
\end{eqnarray}
where $C_{jj'}$ denotes the common copula of $(X_{n,j},X_{n,j'})$, $n\geq 1$. Moreover,
\begin{eqnarray}\label{plambda1.2}
\begin{array}{rl}
1-\frac{1-t}{F_{j'}(c_{j'}^{-r}w_{j't})}\leq 2-\frac{1}{t}\bigg(1-C_{jj'}( 1-t,F_{j'}(c_{j'}^{-r}w_{j't}))\frac{1-t}{F_{j'}(c_{j'}^{-r}w_{j't})}\bigg)\leq 2-\frac{1}{t}\bigg(1-\frac{(1-t)^2}{F_{j'}(c_{j'}^{-r}w_{j't})}\bigg).
\end{array}
\end{eqnarray}
\end{pro}
\dem We have that
\begin{eqnarray}\nn%
\begin{array}{rl}
&\dst\lim_{t\downarrow 0}\frac{P( F_j(X_{1,j}) > 1-t,F_{j'}(X_{1+r,j'}) > 1-t)}
{P( F_j(X_{1,j}) > 1-t)}\vspace{0.35cm}\\
=& 2-\dst\lim_{t\downarrow 0}\frac{1}{t}\big(1-P( F_j(X_{1,j}) \leq 1-t,F_{j'}(X_{1+r,j'}) \leq 1-t)\big)
\vspace{0.35cm}\\
=& 2-\dst\lim_{t\downarrow 0}\frac{1}{t}\bigg(1-P( F_j(X_{1,j}) \leq 1-t,F_{j'}(X_{1,j'}) \leq F_{j'}(c_{j'}^{-r}w_{j't}))\prod_{i=1}^rG_{j'}(w_{j't}/c_{j'}^{r-i})\bigg)
\end{array}
\end{eqnarray}
Now the first result  follows from (\ref{fmar}). The second assertion is a consequence of the Fr\'{e}chet-Hoeffding copula bounds, i.e., $\max(u_1+u_2-1,0)\leq C(u_1,u_2)\leq \min(u_1,u_2)$, for all $(u_1,u_2)\in [0,1]^2$. \fdem\\\\

In the following we will state some consequences of this result, considering different situations for the domains of attraction of $F_{j'}$.

\begin{pro}\label{plambda2}
Under the conditions of Proposition \ref{plambda1}, we have
\begin{eqnarray}\nn
\begin{array}{rl}
\lambda_{jj'}^{(r)}(\mb{X})=
2-\dst\lim_{t\downarrow 0}\frac{1}{t}\bigg(1-C_{jj'}\big(1-t,1-t{c_{j'}^{r\alpha_{j'}}}\big)
\frac{1-t}{1-t{c_{j'}^{r\alpha_{j'}}}}\bigg)
\end{array}
\end{eqnarray}
if $F_{j'}\in\mathcal{D}(\Phi_{\alpha_{j'}})$. Moreover, $0\leq \lambda_{jj'}^{(r)}(\mb{X})\leq c_{j'}^{\alpha_{j'}r}$.
\end{pro}
\dem The result is straightforward since we have $F_{j'}(c_{j'}^{-r}w_{j't})\dst\mathop{\sim}_{t\downarrow 0}1-t{c_{j'}^{r\alpha_{j'}}}$ (see, for instance, Proposition 3.3 in Ferreira and Canto e Castro \cite{mf+lcc} 2008). \fdem\\

\begin{pro}\label{plambda3}
Under the conditions of Proposition \ref{plambda1}, we have $\lambda_{jj'}^{(r)}(\mb{X})=0$ and $\eta_{jj'}^{(r)}(\mb{X})=1/2$, whenever $F_{j'}$ has positive finite right end-point, with $r\in\en$.
\end{pro}
\dem  Observe that $F_{j'}(F_{j'}^{-1}(1-t)c_{j'}^{-r})\dst\mathop{\sim}_{t\downarrow 0}1$, and thus
\begin{eqnarray}\nn
\begin{array}{rl}
2-\frac{1}{t}\big(1-C_{jj'}(1-t,1)(1-t)\big){\dst\mathop{\sim}_{t\downarrow 0} } 2-\frac{1}{t}(1-(1-t)^2){\dst\mathop{\sim}_{t\downarrow 0} } t. \,\,\,\,\square
\end{array}
\end{eqnarray}\\

\begin{pro}\label{plambda4}
Under the conditions of Proposition \ref{plambda1}, we have $\lambda_{jj'}^{(r)}(\mb{X})=0$ and $1/2\leq \eta_{jj'}^{(r)}(\mb{X})\leq \max(1/2,c_{j'}^{rk})$ whenever $F_{j'}(c_{j'}^{-r}w_{j't})\sim 1-t^{c_{j'}^{-rk}}$, for $k>0$ and $r\in\en$.
\end{pro}
\dem Just observe that the left and right hand-side of (\ref{plambda1.2}) approximates, respectively, $t$ and $t+t^{c_{j'}^{-rk}-1}$, as $t\downarrow 0$. \fdem \\\\

Examples of d.f.'s satisfying $F_{j'}(c_{j'}^{-r}w_{j't})\sim 1-t^{c_{j'}^{-rk}}$, $k>0$, include, e.g., Weibull of minimums (with d.f. $F(x)=1-\exp(-x^k)$) and Exponential ($k=1$).


In a max-autoregressive context, the non negative associated tail independence ($1/2\leq \eta<1$)  can also be described through a $p$ARMAX process (Ferreira and Canto e Castro \cite{mf+lcc,lccmf2}, 2008/10), whose logarithm corresponds to ARMAX.\\

\begin{figure}
\begin{center}
\includegraphics[width=12cm,height=4cm]{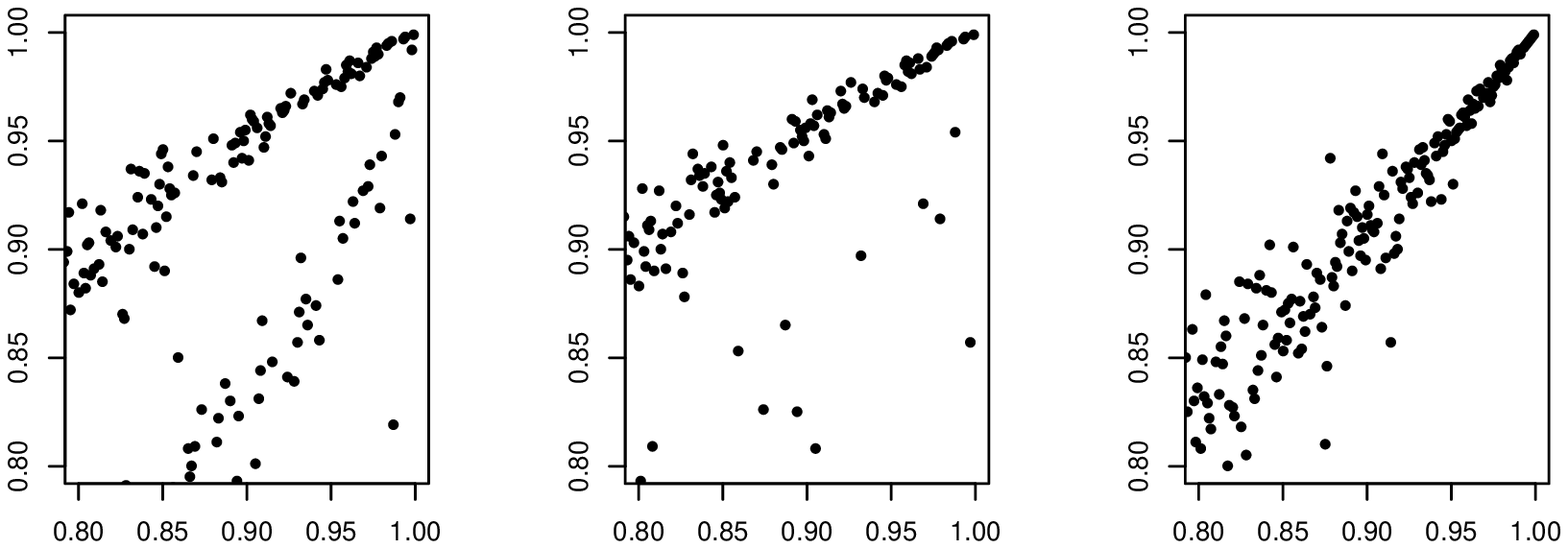}
\includegraphics[width=12cm,height=4cm]{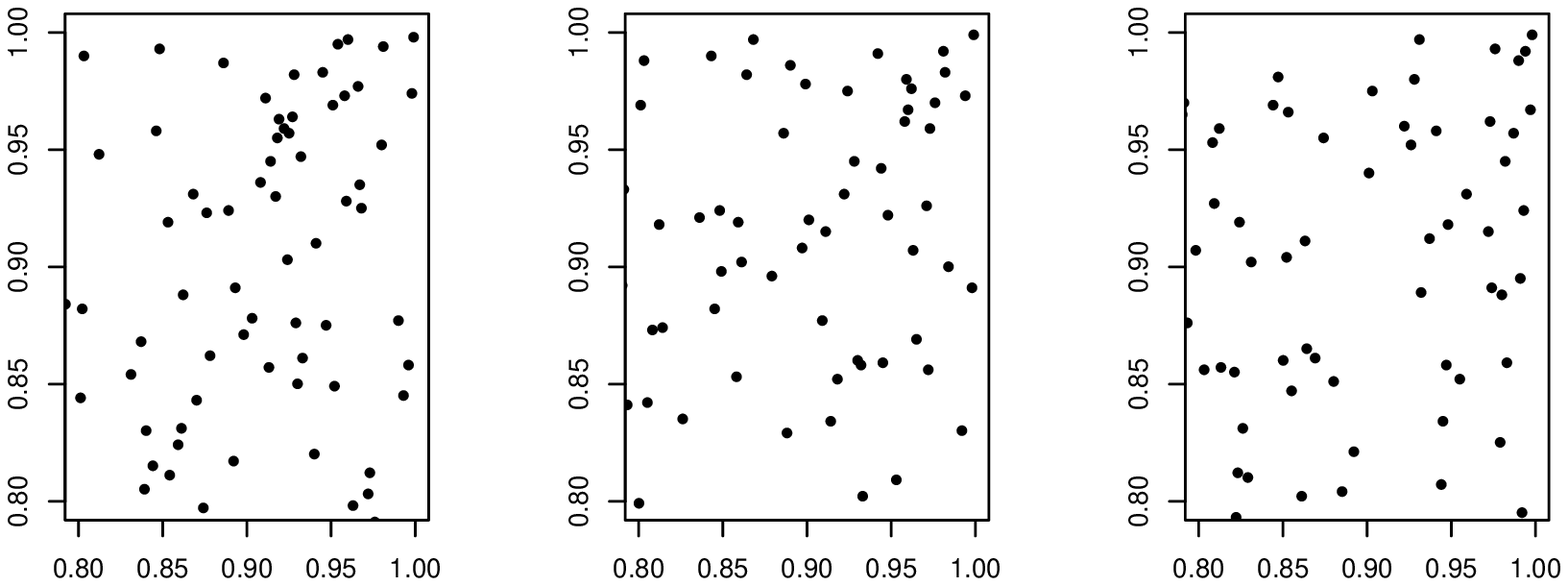}
\end{center}
\caption{Scatter-plots for model $(X_1,X_2,X_3)$ based on Gumbel's copula and marginals ARMAX, respectively, $c=0.8,0.1,0.1$ and  innovations distributed as  unit Fr\'{e}chet,  Exponential and  Uniform; Left to right: points of $(X_1,X_2)$, $(X_1,X_3)$ and $(X_2,X_3)$ with Gumbel's copula dependence parameter $\gamma=0.1$ (strong dependence) on the top and Gumbel's copula dependence parameter $\gamma=0.9$ (weak dependence) on the bottom. \label{fig1}}
\end{figure}

\begin{figure}
\begin{center}
\includegraphics[width=12cm,height=11cm]{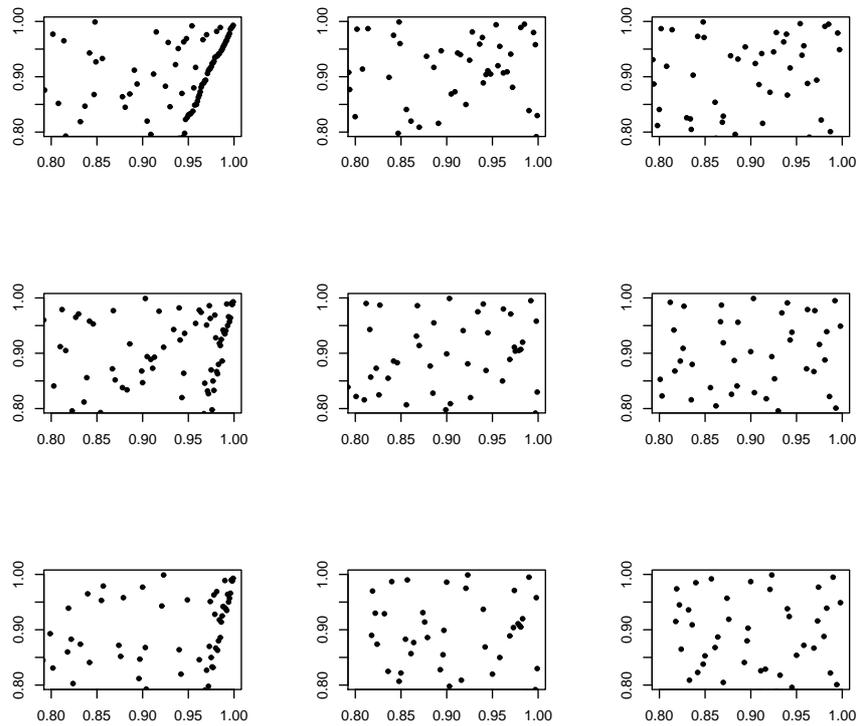}
\end{center}
\caption{Cross-sectional scatter-plots for model $(X_1,X_2,X_3)$ based on Gumbel's copula and marginals ARMAX($0.5$) with innovations distributed as, respectively, unit Fr\'{e}chet,  Exponential and Uniform; Left to right and top to bottom: points of $(X_1,X_1^{(2)})$, $(X_1,X_2^{(2)})$,  $(X_1,X_3^{(2)})$,  $(X_2,X_1^{(2)})$, $(X_2,X_2^{(2)})$, $(X_2,X_3^{(2)})$, $(X_3,X_1^{(2)})$, $(X_3,X_2^{(2)})$ e $(X_3,X_3^{(2)})$, where $X_j^{(r)}$ denotes the $j^{th}$ marginal lag-$r$ apart. \label{fig2}}
\end{figure}

An illustration of the tail dependence between the marginals of $\{\mb{X}_n\}_{n\geq 1}$, for the three cases of domains of attraction can be seen in Figure \ref{fig1}. Observe that the copula's dependence is determinant: dependence is evident whenever a strong dependent copula is used, whilst a  weak dependent copula leads to an almost random  scatter-plot. Figure \ref{fig2} illustrates cross-sectional tail dependence  of $\{\mb{X}_n\}_{n\geq 1}$, considering again the three domains. Observe the presence of some tail dependence for random pairs $(X_j,X_{j'}^{(r)})$ whenever the lag-$r$ apart $j'^{th}$ marginal is Fréchet (first column plots) corroborating Proposition \ref{plambda2}. An almost randomness can be seen in the other scatter-plots which is consistent with Propositions \ref{plambda3} and \ref{plambda4}.

%
%
%

\section{Marginal parameters estimation}\label{sestim}

In this section we shall focus on the marginal ARMAX autoregressive parameter and the marginal tail index.

The following result allow us to state an estimator for the ARMAX parameter $c_j$ , $j\in D$.
\begin{pro}\label{pcj}
If $F_{0,j}$ and $G_j$ are unit Fréchet d.f.'s, $j\in D$, then \begin{eqnarray}\label{cj}
\begin{array}{rl}
\dst c_j=2-\frac{1}{E\big(e^{-X_j^{-1}}\big)}
\end{array}
\end{eqnarray}
where $X_j$ is a r.v. with the stationary distribution of $\{X_{n,j}\}_{n\geq 1}$.
\end{pro}
\dem We shall use the result of Proposition 3.1 in Ferreira and Ferreira (\cite{hf+mf3}, 2012). More precisely, considering $F(x)=e^{-x^{-1}}$ and $s\in\en$, we have
\begin{eqnarray}\label{pcj1}
\begin{array}{rl}
E(F(X_{n,j})^s)=&E(F(c_j^nX_{0,j}\vee \bigvee_{i=1}^n c_j^{n-i}Y_{i,j})^s)\vspace{0.35cm}\\
=&\dst\frac{-\log F_{(X_{0,j},Y_{1,j},...,Y_{n,j})}(c_j^{-n},c_j^{-n+1},...,c_j^{-1},1)}{s-\log F_{(X_{0,j},Y_{1,j},...,Y_{n,j})}(c_j^{-n},c_j^{-n+1},...,c_j^{-1},1)}\vspace{0.35cm}\\
=&\dst\frac{\sum_{i=0}^nc_j^{n-i}}{s+\sum_{i=0}^nc_j^{n-i}}.
\end{array}
\end{eqnarray}
Assuming that $\{X_{n,j}\}_{n\geq 1}$ is stationary and taking limits in both of the members with $s=1$, we have
\begin{eqnarray}\nn
\begin{array}{rl}
E\big(e^{-X_j^{-1}}\big)=\dst\frac{\frac{1}{1-c_j}}{1+\frac{1}{1-c_j}},
\end{array}
\end{eqnarray}
leading to the assertion. \fdem\\

As a consequence of this result, we verify that if $\{X_{n,j}\}_{n\geq 1}$ is stationary than $E\big(e^{-X_{j}^{-1}}\big)\in (1/2,1)$.\\

From (\ref{cj}) we derive the estimator
\begin{eqnarray}\nn
\begin{array}{rl}
\dst \widehat{c_j}=2-\frac{1}{\overline{U_j}}
\end{array}
\end{eqnarray}
where $\overline{U_j}=\frac{1}{n}\sum_{i=1}^n e^{-X_{i,j}^{-1}}$. Observe that no definite result can be obtained for $\overline{U_j}\leq 1/2$, which may be an indication of an unsuitable model's choice.\\

Based on (\ref{mar1}), we have that $X_{n,j}=\bigvee_{i=1}^\infty c_j^{i} Y_{n-i,j}$ is the unique stationary solution of recursion (\ref{armax}) (Davis and Resnick \cite{dav+res}, 1989; Proposition 2.2). Therefore, an ARMAX process is ergodic (Stout \cite{stout74}, 1974; Theorem 3.5.8) and, since $E(|e^{-X_j^{-1}}|) < \infty$, we have $\overline{U_j}\to E\big(e^{-X_j^{-1}}\big)$ almost surely (see, e.g. ergodic theory in Billingsley \cite{bil} 1995). Thus estimator $\widehat{c_j}$ is strongly consistent. The asymptotic normality is stated in the next result.

\begin{pro}\label{pecj}
Under the conditions of Proposition \ref{pcj}, we have $\sqrt{n}(\overline{U_j}-e^{-X_{j}^{-1}})\to N(0,\sigma^2)$
where
\begin{eqnarray}\nn
\begin{array}{rl}
\dst {\sigma^2}=\frac{1}{3-2c_j}-\Big(\frac{1}{2-c_j}\Big)^2+2\sum_{r=1}^\infty \Big(\frac{1-c_j^r}{(2-c_j)(2-c_j-c_j^r-c_j^{r+1})}-\Big(\frac{1}{2-c_j}\Big)^2\Big).
\end{array}
\end{eqnarray}
Moreover $\sqrt{n}(\widehat{c_j}-c_j)\to N\big(0,\sigma^2(3-2c_j)\big)$.
\end{pro}
\dem
The asymptotic normality also holds given the strong-mixing dependence structure with variance given by (Billingsley, \cite{bil}, 1995)
\begin{eqnarray}\nn
\begin{array}{rl}
\dst {\sigma^2}=var\big(e^{-X_j^{-1}}\big)+2\sum_{r=1}^\infty cov\big(\big(e^{-X_j^{-1}}\big)\big(e^{-X_{j+r}^{-1}}\big)\big).
\end{array}
\end{eqnarray}
Since $G_j$ and $F_{0,j}$ are unit Fréchet, according to the stationarity relation in (\ref{fmar}), then
\begin{eqnarray}\nn
\begin{array}{rl}
F_j(x)=e^{-x^{-1}\frac{1}{1-c_j}},\,x>0,\,j\in D.
\end{array}
\end{eqnarray}
For each $r\in \en$, the joint d.f. of $(X_{n,j},X_{n+r,j})$ is given by
\begin{eqnarray}\nn
\begin{array}{rl}
F_{(X_{n,j},X_{n+r,j})}(x,y)
=P(X_{n,j}\leq x\wedge yc_j^{-r})P(Y_{n+1,j}\leq yc_j^{-r+1},...,Y_{n+r,j}\leq yc_j),\,x>0,\,j\in D.
\end{array}
\end{eqnarray}
whose joint density, for $y>xc_j^{-r}>0 $, is derived as
\begin{eqnarray}\nn
\begin{array}{rl}
\dst\frac{\partial}{\partial x}\frac{\partial}{\partial y}F_{(X_{n,j},X_{n+r,j})}(x,y)=&\dst\frac{\partial}{\partial x}\frac{\partial}{\partial y}\bigg(e^{-x^{-1}\frac{1}{1-c_j}}\prod_{i=1}^r e^{-y^{-1}{c_j^{r-i}}}\bigg)\vspace{0.35cm}\\
=&\dst\frac{1}{x^2(1-c_j)}e^{-x^{-1}\frac{1}{1-c_j}}
\frac{1-c_j^r}{y^2(1-c_j)}e^{-y^{-1}\frac{1-c_j^r}{1-c_j}},\,j\in D.
\end{array}
\end{eqnarray}
Therefore, and after some calculations, we obtain
\begin{eqnarray}\nn
\begin{array}{rl}
E\big(e^{-X_j^{-1}}e^{-X_{j+r}^{-1}}\big)=&\dst\int_0^\infty\int_0^{yc_j^{-r}}e^{-x^{-1}}e^{-y^{-1}}\frac{\partial}{\partial x}\frac{\partial}{\partial y}F_{(X_{n,j},X_{n+r,j})}(x,y)dxdy\vspace{0.35cm}\\
=&\dst\frac{1-c_j^r}{(2-c_j)(2-c_j-c_j^r-c_j^{r+1})}.
\end{array}
\end{eqnarray}
Now just observe that, by (\ref{pcj1}),
\begin{eqnarray}\nn
\begin{array}{rl}
E\Big(\big(e^{-X_j^{-1}}\big)^2\Big)=\dst\frac{\frac{1}{1-c_j}}{2+\frac{1}{1-c_j}}=\frac{1}{3-2c_j}.
\end{array}
\end{eqnarray}
Considering $g(x)=2-1/x$, we have $\Big[g'\Big(E\big(e^{-X_j^{-1}}\big)\Big)\Big]^2=E\big(e^{-X_j^{-1}}\big)^{-2}$
and, by the Delta Method, the second assertion holds. \fdem\\


Other estimators can be found in literature. A strongly consistent  estimator was earlier proposed in Davis and Resnick (\cite{dav+res}, 1989):
\begin{eqnarray}\nn
\begin{array}{rl}
\dst \widetilde{c}_j^*=\bigwedge_{i=2}^n\frac{X_i}{X_{i-1}}.
\end{array}
\end{eqnarray}
Another estimator, with a quite similar expression to our proposal, was derived in Lebedev (\cite{leb}, 2008). More precisely, for unit Fréchet marginals $F_j$ and $F_{0,j}(x)=G_j(x)=e^{-x^{-1}\frac{1}{1-c_j}}$, then
\begin{eqnarray}\nn
\begin{array}{rl}
\dst c_j=2-\frac{1}{p_j},
\end{array}
\end{eqnarray}
with $p_j=P(X_{2,j}\leq X_{1,j}) \in (1/2,1)$, and thus
\begin{eqnarray}\nn
\begin{array}{rl}
\dst \widetilde{c}_j=2-\frac{1}{\widetilde{p_j}}
\end{array}
\end{eqnarray}
where $\widetilde{p_j}=(n-1)^{-1}\sum_{i=1}^{n-1}\mathds{1}_{\{X_{i+1,j}\leq X_{i,j}\}}$. Note that a similar restriction to our method must be considered, i.e., $1/2<\widetilde{p_1}<1$. The consistency and asymptotic normality of this estimator can be seen in Ferreira (\cite{mf12}, 2012).\\


In what concerns the tail index $\alpha_j$ of each marginal $j\in D$, it can be estimated through tail index estimators already stated in literature as Hill (in case $\alpha_j>0$), Pickand's, maximum likelihood, moments or generalized moments estimator, whose asymptotic properties of consistency and normality still hold under an ARMAX dependence structure (Ferreira and Canto e Castro \cite{mf+lcc} 2008, Proposition 3.4).

\end{document}